\documentclass[journal,twoside,web]{ieeecolor}
\usepackage{generic}
\usepackage{cite}
\usepackage{amsmath,amssymb,amsfonts}
\usepackage{algorithmic}
\usepackage{graphicx}
\usepackage{textcomp}
\usepackage{multirow}

\newtheorem{thm}{Theorem}

\def\BibTeX{{\rm B\kern-.05em{\sc i\kern-.025em b}\kern-.08em
    T\kern-.1667em\lower.7ex\hbox{E}\kern-.125emX}}
\begin{document}
\title{Embedded Point Iteration Based Recursive Algorithm for Online Identification of Nonlinear Regression Models}
\author{Guang-Yong Chen,
        Min Gan,~\IEEEmembership{~Senior Member,~IEEE},
        Jing Chen,
        Long Chen\\
\thanks{G.-Y. Chen \& M. Gan are with the College of Computer and Data Science, Fuzhou University, Fuzhou 350116, China
(e-mail: cgykeda@mail.ustc.edu.cn; aganmin@aliyun.com).}
\thanks{J. Chen is with School of Science, Jiangnan University, Wuxi, 214122, China (e-mail: chenjing1981929@126.com)}
\thanks{L. Chen is with the Faculty of Science and Technology, University of Macau, Macau 99999, China (e-mail: longchen@umac.mo) }
}

\maketitle

\begin{abstract}
This paper presents a novel online identification algorithm for nonlinear regression models. The online identification problem is challenging due to the presence of nonlinear structure in the models.
Previous works usually ignore the special structure of nonlinear regression models, in which the parameters can be partitioned into a linear part and a nonlinear part. In this paper, we develop an efficient recursive algorithm for nonlinear regression models based on analyzing the equivalent form of variable projection (VP) algorithm. By introducing the embedded point iteration (EPI) step, the proposed recursive algorithm can properly exploit the coupling relationship of linear parameters and nonlinear parameters. In addition, we theoretically prove that the proposed algorithm is mean-square bounded. Numerical experiments on synthetic data and real-world time series verify the high efficiency and robustness of the proposed algorithm.
\end{abstract}

\begin{IEEEkeywords}
Nonlinear regression models, online identification, parameter estimation, variable projection.
\end{IEEEkeywords}

\section{Introduction}
\label{sec:introduction}
\IEEEPARstart{T}{ime} series analysis plays a vital important role in scientific research and engineering applications. At early stage, the linear models were used as convenient tools for time series modeling and forecasting, however, they usually fail in capturing some important nonlinear behaviors presented in the time series, e.g., time irreversibility, asymmetry, and self-sustained stochastic cyclical behavior. Aware of these, researchers focused on the nonlinear models and proposed various models in different fields. For example, %Tong \cite{tong2012threshold} proposed the threshold autoregressive (TAR) model that has been widely used in neuroscience, hydrology and finance.
Ozaki \& Oda \cite{ozaki1977non} used the exponential autoregressive (Exp-AR) model to model the ship rolling data.
%, which successfully reveals the jump phenomena, amplitude dependent frequency shifts and perturbed limit cycles presented in the data.
In \cite{priestley1980state}, Priestly constructed a general class of nonlinear models, named state-dependent autoregressive (SD-AR) model, to forecast and indicate specific types of nonlinear behavior. With simple topological structure and strong learning capability, the radial basis function (RBF) neural network offers a viable alternative to the traditional models. Vesin \cite{vesin1993amplitude} and Shi \emph{et al.} \cite{shi1999nonlinear} adopted an RBF approach to SD-AR modeling and developed the RBF-AR model. Peng \emph{et al.} \cite{peng2003tnn} further extended the idea of RBF-AR to the systems that have several external inputs, and proposed the RBF-ARX model.
%which has been employed in modeling the magnetic levitation system \cite{qin2014modeling}, nitrogen oxide ($\text{NO}_x$) decomposition process \cite{peng2004cep}, thermal power plant \cite{peng2011amm}, and etc.

Most of these models mentioned above fall into a class of nonlinear regression models that take the form:
\begin{equation}\label{1-1}
  f(\boldsymbol a,\boldsymbol c; \boldsymbol x)=\sum\limits_{i=1}^nc_i\phi_{i}(\boldsymbol a;\boldsymbol x),
\end{equation}
where $\boldsymbol a=(a_1,\cdots,a_k)^\text T$, $\boldsymbol c=(c_1, c_2, \cdots, c_n)^\text T$ are the parameters to be estimated, $\boldsymbol x$ is the state vector that may consist of the input or output, or any other explanatory variables in the system, and $\phi_i(\cdot)$ is a differentiable nonlinear  function.

For offline identification, various approaches presented in the literature can be used for the parameter estimation of the nonlinear regression models, e.g., the joint optimization algorithm \cite{okatani2007wiberg}, alternating optimization method \cite{okatani2007wiberg,okatani2011efficient}, variable projection (VP) algorithm \cite{golub1973differentiation,golub2003separable}, inner point solver \cite{aravkin2013sparse}, and structured nonlinear parameter optimization method (SNPOM) \cite{zeng2017regularized}. The VP algorithm, which takes advantage of the separable structure of nonlinear regression models, is more powerful and usually yields a better-conditioned problem \cite{sjoberg1997separable,chen2021insights}. Empirical evidences \cite{chen2019regularized,aravkin2012estimating,gan2017some,erichson2020sparse}
have proven that the VP algorithm achieves faster convergence rate and is less sensitive to the initial value of the chosen parameters.

But in real applications, the data samples usually arrive sequentially, or the system is time-varying. The VP algorithm for batch data in offline scenario could not be applied to the online system. Most of the recursive algorithms presented in the literature \cite{ngia2000efficient,asirvadam2002separable,shamsudin2014recursive,chen2020modified} ignore the special structure of nonlinear regression models. For example,  Ngia \& Sjoberg \cite{ngia2000efficient} treated all the parameters of the neural networks the same and proposed the recursive Gauss-Newton (RGN) and Levenberg-Marquardt (RLM) algorithms for the online learning. Asivadam \emph{et al.} \cite{asirvadam2002separable} presented a hybrid recursive Levenberg-Marquardt for feedforward neural networks, which is a strategy that optimizes linear and nonlinear parameters alternately.

Recently, Gan \emph{et al.} \cite{gan2019adaptive,gan2020recursive} incorporated the VP step into the online identification algorithm for separable nonlinear optimization problems. However, they did not properly resolve the coupling relationship between linear parameters and nonlinear parameters. In \cite{gan2019adaptive}, the authors just used a few data samples to estimate the nonlinear parameters, and then in the subsequent learning process, only the linear parameters are updated. In \cite{gan2020recursive}, the VP step was introduced to update the nonlinear parameters, but the authors employed the least-squares method to eliminate the linear parameters only using the previous $p$ data samples (from time $t-p+1$ to $t$). This approach has two main defects: 1) If $p$ is small, the linear parameter is only determined by a small number of observations, which greatly limits the identification accuracy of the algorithm; 2) if $p$ is large, calculating generalized inverse consumes much computational cost. Both problems are contrary to the design principle of online learning algorithm.

%The classical implementation of VP algorithm for batch data samples is not suitable for recursive learning algorithms for separable nonlinear models.
For online learning, the most difficult problem in applying the VP strategy is how to solve the coupling relationship between the linear and nonlinear parameters in the recursive process. This motivates us to break through the framework of classic VP algorithm proposed by Golub \& Pereyra \cite{golub1973differentiation}, and explore new ideas for designing online VP algorithm. The main contribution of this paper is to analyze VP algorithm from another perspective and then to employ its equivalent form---the embedded point iteration (EPI) algorithm to design a recursive learning algorithm for online learning of the nonlinear regression models. The proposed recursive algorithm considers the coupling relationship of the linear parameters and nonlinear parameters during the recursive process. In addition, we theoretically prove that the proposed recursive algorithm is mean-square bounded. It is conceivable that the proposed recursive algorithm will be of great important because of the widespread of nonlinear regression models in the real-world applications.

The paper proceeds as follows. The equivalent form of VP algorithm is analyzed in Section \uppercase\expandafter{\romannumeral2}.  In Section \uppercase\expandafter{\romannumeral3}, we introduce the recursive Gauss-Newton (RGN) algorithm for online learning of nonlinear regression models. The proposed recursive algorithm based on the EPI method is given in Section \uppercase\expandafter{\romannumeral4}. Numerical experiments are carried out to confirm the performance of the proposed recursive algorithm in Section \uppercase\expandafter{\romannumeral5}. Finally, we draw the conclusion of this paper.

\section{Equivalent form of VP algorithm}
In this paper, we consider the nonlinear regression model (\ref{1-1}). Given observations $\{y_1,\cdots,y_m\}$, the parameter identification of model (\ref{1-1}) can be formulized as the minimization problem
\begin{equation}\label{m2-1}
  \min\limits_{\boldsymbol a, \boldsymbol c}\  r(\boldsymbol a, \boldsymbol c)=\sum\limits_{i=1}^m(y_i-\sum\limits_{j=1}^nc_j\phi_j(\boldsymbol a;\boldsymbol x_i))^2,
\end{equation}
which can be written in a matrix form
\begin{equation}\label{snlls1}
  \min\limits_{\boldsymbol a, \boldsymbol c}\  r(\boldsymbol a, \boldsymbol c)=\left\|\mathbf \Phi(\boldsymbol a)\boldsymbol c-\boldsymbol y\right\|_2^2,
\end{equation}
where $\boldsymbol y=(y_1,\cdots,y_m)^\text T$ is the observation vector, the variables $\boldsymbol c=(c_1,\cdots,c_n)^\text T$ and $\boldsymbol a$ are regarded as linear/nonlinear parameters, and $\left(\mathbf \Phi(\boldsymbol a)\right)_{ij}=\phi_j(\boldsymbol a;\boldsymbol x_i)$. Golub \& Pereyra \cite{golub1973differentiation} referred the data fitting problem with model (\ref{1-1}) as separable nonlinear least squares (SNLLS) problems. Taking advantage of the special structure of SNLLS problems, they proposed an efficient VP strategy, which considers the coupling relationships of the linear parameters and nonlinear parameters. The basic idea of VP algorithm is to eliminate the linear parameters, and then to minimize the reduced function that only contains the nonlinear parameters.

For fixed parameter $\boldsymbol a$, the optimization problem (\ref{snlls1}) is convex. The optimal solution of linear parameters can be obtained by solving the linear least squares problem:
\begin{equation}\label{E-2}
  \hat{\boldsymbol c}=\arg\min\limits_{\boldsymbol c}\ r(\boldsymbol a, \boldsymbol c)=\mathbf\Phi^\dagger\boldsymbol y,
\end{equation}
where $\mathbf\Phi^\dagger=(\mathbf\Phi^\text T\mathbf\Phi)^{-1}\mathbf\Phi^\text T$ is the Moore-Penrose inverse of $\mathbf\Phi$. Replacing $\boldsymbol c$ in (\ref{snlls1}) with (\ref{E-2}) yields a reduced function:
\begin{equation}\label{E-3}
  r_2(\boldsymbol a)=\left\|\boldsymbol r_2\right\|_2^2=\left\|\mathbf P_{\mathbf\Phi}^\bot\boldsymbol y\right\|_2^2=\left\|(\mathbf I-\mathbf \Phi\mathbf\Phi^\dagger)\boldsymbol y\right\|_2^2,
\end{equation}
where $\mathbf P_{\mathbf\Phi}^\bot=\mathbf I-\mathbf \Phi\mathbf\Phi^\dagger$ is a project operator.
%that projects a vector $\boldsymbol y$ to the orthogonal complement of the column space of $\mathbf\Phi$.
Compared with the problem (\ref{snlls1}), the reduced objective function is more complicated, although it contains fewer parameters to be estimated. Fortunately, Golub \& Pereyra \cite{golub1973differentiation} proved that the original problem and the reduced problem achieve the same global minimizer and obtained the Jacobian matrix of the reduced function
\begin{equation}\label{jac1}
  \mathbf J_{\text{GP}}=-\mathbf P_{\boldsymbol\Phi}^\bot\text D\boldsymbol\Phi\boldsymbol\Phi^-\boldsymbol y-(\mathbf P_{\boldsymbol\Phi}^\bot\text D\boldsymbol\Phi\boldsymbol\Phi^-)^\text T\boldsymbol y,
\end{equation}
where $\text D \boldsymbol\Phi$ represents the Fr\'{e}chet derivative of $\boldsymbol\Phi$, and $\boldsymbol\Phi^-$ is generalized inverse of $\boldsymbol\Phi$, which satisfies $(\boldsymbol\Phi\boldsymbol\Phi^-)^\text T=\boldsymbol\Phi\boldsymbol\Phi^-$ and $\boldsymbol\Phi\boldsymbol\Phi^-\boldsymbol\Phi=\boldsymbol\Phi$. %Note that the Moore-Penrose inverse in (\ref{E-3}) is a special generalized inverse.
In \cite{kaufman1975variable}, Kaufman proposed a simplified Jacobian matrix, which drops the second term in (\ref{jac1})
\begin{equation}\label{jac2}
  \mathbf J_{\text{Kau}}=-\mathbf P_{\boldsymbol\Phi}^\bot\text D\boldsymbol\Phi\boldsymbol\Phi^-\boldsymbol y.
\end{equation}

Empirical evidence \cite{ruhe1980algorithms,golub2003separable,gan2017some} suggests that Kaufman's simplified VP algorithm has similar performance  to the full form.
%Without calculating the second term in (\ref{jac1}), Kauman's simplified VP algorithm saves the computational cost.
%Due to the high efficiency of VP algorithm, a variety of applications from system identification, image restoration, electrical engineering, environmental sciences have been reported in the literature.
The complexity of the reduced function (\ref{E-3}) and Jacobian matrices make it difficult to apply the classical VP algorithm to the online identification of nonlinear regression models directly. To overcome this difficulty, we examine the VP algorithm from another perspective in the following part.

Denote $\mathbf J_{\boldsymbol a}=\frac{\partial\boldsymbol r}{\partial\boldsymbol a}$, $\mathbf J_{\boldsymbol c}=\frac{\partial\boldsymbol r}{\partial\boldsymbol c}$, and $\boldsymbol r$ as the residual of the minimization problem (\ref{snlls1}). The update direction $(\Delta \boldsymbol a, \Delta \boldsymbol c)$ can be derived by solving
\begin{equation}\label{E-4}
[\mathbf J_{\boldsymbol a}, \mathbf J_{\boldsymbol c}]\left[
                                                                       \begin{array}{c}
                                                                         \Delta \boldsymbol a \\
                                                                         \Delta \boldsymbol c \\
                                                                       \end{array}
\right]=-\boldsymbol r.
\end{equation}
Applying the linear least squares method to (\ref{E-4}), we have
\begin{small}
\begin{equation}\label{E-5}
  \left[
          \begin{array}{cc}
            \mathbf J_{\boldsymbol a}^{\text{T}}\mathbf J_{\boldsymbol a}\  &\mathbf J_{\boldsymbol a}^{\text{T}}\mathbf J_{\boldsymbol c}  \\
            \\
            \mathbf J_{\boldsymbol c}^{\text{T}}\mathbf J_{\boldsymbol a}\  & \mathbf J_{\boldsymbol c}^{\text{T}}\mathbf J_{\boldsymbol c} \\
          \end{array}
  \right]
   \left[
                                                                       \begin{array}{c}
                                                                         \Delta \boldsymbol a \\
                                                                         \Delta \boldsymbol c \\
                                                                       \end{array}
\right]=
\left[
                                                                       \begin{array}{c}
                                                                         -\mathbf J_{\boldsymbol a}^{\text{T}}\boldsymbol r \\
                                                                         \\
                                                                         -\mathbf J_{\boldsymbol c}^{\text{T}}\boldsymbol r \\
                                                                       \end{array}
\right].
\end{equation}
\end{small}

The method of calculating the update direction by solving (\ref{E-5}) is referred to as the joint optimization algorithm. For fixed parameters $\boldsymbol a$, the optimal value of $\boldsymbol c$ is the solution of linear least squares problem (\ref{E-2}), which yields $\frac{\partial r}{\partial \boldsymbol c}=\mathbf J_{\boldsymbol c}^{\text{T}}\boldsymbol r=0$. Inserting it into (\ref{E-5}), we have
\begin{small}
\begin{equation}\label{E-6}
  \left[
          \begin{array}{cc}
            \mathbf J_{\boldsymbol a}^{\text{T}}\mathbf J_{\boldsymbol a}\  &\mathbf J_{\boldsymbol a}^{\text{T}}\mathbf J_{\boldsymbol c}  \\
            \\
            \mathbf J_{\boldsymbol c}^{\text{T}}\mathbf J_{\boldsymbol a}\  & \mathbf J_{\boldsymbol c}^{\text{T}}\mathbf J_{\boldsymbol c} \\
          \end{array}
  \right]
   \left[
                                                                       \begin{array}{c}
                                                                         \Delta \boldsymbol a \\
                                                                         \Delta \boldsymbol c \\
                                                                       \end{array}
\right]=
\left[
                                                                       \begin{array}{c}
                                                                         -\mathbf J_{\boldsymbol a}^{\text{T}}\boldsymbol r \\
                                                                         \\
                                                                         0\\
                                                                       \end{array}
\right].
\end{equation}
\end{small}
This algorithm is named as embedded point iteration (EPI), which has been employed in bundle adjustment \cite{hyeong2017revisiting}.
%Using the Schur complement of a matrix, we can obtain the update direction
%\begin{eqnarray}\label{E-7}
%  \Delta \boldsymbol a&=&-(\mathbf J_{\boldsymbol a}^\text{T}(\mathbf I_m-\mathbf J_{\boldsymbol c}\mathbf J_{\boldsymbol c}^\dagger)\mathbf J_{\boldsymbol a})^{-1}\mathbf J_{\boldsymbol a}^\text{T}\boldsymbol r\nonumber\\
%  &=&
%\end{eqnarray}
The equation (\ref{E-6}) yields
\begin{equation}\label{E-7}
  \mathbf J_{\boldsymbol c}^\text{T}\mathbf J_{\boldsymbol c}\Delta \boldsymbol c+\mathbf J_{\boldsymbol c}^\text{T}\mathbf J_{\boldsymbol a} \Delta \boldsymbol a=0,
\end{equation}
that is,
\begin{eqnarray}\label{E-8}
  \frac{\Delta \boldsymbol c}{\Delta \boldsymbol a}=(\mathbf J_{\boldsymbol c}^\text{T}\mathbf J_{\boldsymbol c})^{-1}\mathbf J_{\boldsymbol c}^\text{T}\mathbf J_{\boldsymbol a}
  =(\mathbf \Phi^\text T\mathbf\Phi)^{-1}\mathbf\Phi^\text T \text D \mathbf\Phi\mathbf\Phi^\dagger\boldsymbol y.
\end{eqnarray}
The second equation holds since $\mathbf J_{\boldsymbol c}=\mathbf\Phi$, $\mathbf J_{\boldsymbol a}=\text D\mathbf\Phi\mathbf \Phi^\dagger\boldsymbol y$.
The Jacobian matrix of reduced function (\ref{E-3}) can be calculated as
\begin{equation}\label{E-9}
  \mathbf J=\text D \boldsymbol r_2=\mathbf J_{\boldsymbol a}+\mathbf J_{\boldsymbol c}\frac{d \boldsymbol c}{d \boldsymbol a}.
\end{equation}
Replacing $\frac{d \boldsymbol c}{d \boldsymbol a}$ in (\ref{E-9}) with (\ref{E-8}) yields an approximated Jacobian matrix
\begin{eqnarray}\label{E-10}
  \mathbf J\approx\text D\mathbf\Phi\mathbf \Phi^\dagger\boldsymbol y+\mathbf \Phi(\mathbf \Phi^\text T\mathbf\Phi)^{-1}\mathbf\Phi^\text T \text D \mathbf\Phi\mathbf\Phi^\dagger\boldsymbol y
  =-\mathbf P_{\boldsymbol\Phi}^\bot\text D\boldsymbol\Phi\boldsymbol\Phi^\dagger\boldsymbol y,
\end{eqnarray}
which is the same as Kaufman's form (\ref{jac2}). Thus, the EPI approach is equivalent to the Kaufman's VP algorithm.
It is of great significance to expand the VP strategy to online systems.
\section{Recursive Gauss-Newton algorithm}
%Because of the existence of linear and nonlinear parameters that need to be estimated, it is difficult to design an efficient algorithm for the separable nonlinear models.
The RGN algorithm \cite{shamsudin2014recursive} and RLM algorithm \cite{ngia2000efficient} are effective second-order methods for online identification of nonlinear regression models. In this section, we introduce the RGN algorithm that is useful for designing the proposed recursive method.

Denote $\boldsymbol\theta=(\boldsymbol a, \boldsymbol c)$ as the parameters to be estimated in the problem (\ref{snlls1}), the cost function at time $t$ (denoted as $r_t(\boldsymbol \theta)$) can be written as follows
\begin{equation}\label{3-1}
  r_t(\boldsymbol\theta)=\frac{1}{2}\sum\limits_{i=1}^tv_i^2(\boldsymbol\theta)
\end{equation}
where $v_i(\boldsymbol\theta)=y_i-\boldsymbol\phi(\boldsymbol a^{(i)})^\text T\boldsymbol c^{(i)}$; $\boldsymbol a^{(i)}$ and $\boldsymbol c^{(i)}$ are the estimated values of nonlinear parameters and linear parameters at time $i$,  respectively. We consider the Taylor expansion of $r_t(\boldsymbol \theta)$ at $\boldsymbol\theta^{(t-1)}$
\begin{small}
\begin{eqnarray}\label{3-2}
  &&r_t(\boldsymbol\theta)\approx r_t(\boldsymbol\theta^{(t-1)})+r_t'(\boldsymbol\theta^{(t-1)})^\text T(\boldsymbol\theta-\boldsymbol\theta^{(t-1)})\nonumber\\
  &&+\frac{1}{2}(\boldsymbol\theta-\boldsymbol\theta^{(t-1)})^\text T r_t''(\boldsymbol\theta^{(t-1)})(\boldsymbol\theta-\boldsymbol\theta^{(t-1)})\nonumber\\&&+o(||\boldsymbol\theta-\boldsymbol\theta^{(t-1)}||^2),
\end{eqnarray}
\end{small}
where $o(x)$ is the infinitesimal of $x$. Using the GN method to minimize $r_t(\boldsymbol\theta)$ yields
\begin{equation}\label{3-3}
  \boldsymbol \theta^{(t)}=\boldsymbol\theta^{(t-1)}-r_t''(\boldsymbol\theta^{(t-1)})^{-1}r_t'(\boldsymbol\theta^{(t-1)}).
\end{equation}
Specifically, the gradient and Hessian matrix of $r_t(\boldsymbol\theta)$ can be expressed as
\begin{equation}\label{3-4}
  r_t'(\boldsymbol\theta)=\sum\limits_{i=1}^t\boldsymbol g_i(\boldsymbol\theta)v_i(\boldsymbol\theta)=r_{t-1}'(\boldsymbol\theta)+\boldsymbol g_t(\boldsymbol\theta)v_t(\boldsymbol\theta),
\end{equation}
\begin{equation}\label{3-5}
  r_t''(\boldsymbol \theta)=r_{t-1}''(\boldsymbol \theta)+\boldsymbol g_t(\boldsymbol\theta)\boldsymbol g_t(\boldsymbol\theta)^\text T+v_t''(\boldsymbol\theta)v_t(\boldsymbol\theta),
\end{equation}
where $\boldsymbol g_i(\boldsymbol\theta)=\frac{d v_i(\boldsymbol\theta)}{d\boldsymbol\theta}$. As discussed in \cite{ngia2000efficient}, we drop the small term $v_t''(\boldsymbol\theta)v_t(\boldsymbol\theta)$ to obtain an approximated Hessian matrix
\begin{equation}\label{3-6}
  r_t''(\boldsymbol\theta)\approx r_{t-1}''(\boldsymbol\theta)+\boldsymbol g_t(\boldsymbol\theta)\boldsymbol g_t(\boldsymbol\theta)^\text T.
\end{equation}
Denote $\mathbf H_t$ as the approximated Hessian matrix of $r_t(\boldsymbol\theta)$,  and $\mathbf S_t=\mathbf H_t^{-1}$ as the covariance matrix of the estimated parameters $\boldsymbol \theta^{(t)}$. Then,
%The update procedure can be formulated as
%\begin{equation}\label{3-7}
%  \boldsymbol\theta^{(t)}=\boldsymbol\theta^{(t-1)}-\mathbf H_t^{-1}\boldsymbol g_t(\boldsymbol\theta^{(t-1)})v_t(\boldsymbol\theta^{(t-1)}),
%\end{equation}
%and
\begin{equation}\label{3-8}
  \mathbf H_t=\mathbf H_{t-1}+\boldsymbol g_t(\boldsymbol\theta^{(t-1)})\boldsymbol g_t^\text T(\boldsymbol\theta^{(t-1)}),
\end{equation}
%According to the Sherman-Morrison-Woodbury formula \cite{van1983matrix},
%%\begin{small}
%%\begin{equation}\label{3-9}
%%  (\mathbf A+\mathbf U\mathbf V^\text T)^{-1}=\mathbf A^{-1}-\mathbf A^{-1}\mathbf U(\mathbf I+\mathbf V^\text T\mathbf A^{-1}\mathbf U)^{-1}\mathbf V^\text T\mathbf A^{-1},
%%\end{equation}
%%\end{small}
%we obtain the RGN formula
and the RGN formula can be expressed as
\begin{equation}\label{3-10}
  \boldsymbol\theta^{(t)}=\boldsymbol\theta^{(t-1)}-\mathbf S_t\boldsymbol g_t(\boldsymbol\theta^{(t-1)})v_t(\boldsymbol\theta^{(t-1)}),
\end{equation}
\begin{small}
\begin{equation}\label{3-11}
  \mathbf S_t=\mathbf S_{t-1}-\frac{1}{\alpha}\mathbf S_{t-1}\boldsymbol g_t(\boldsymbol\theta^{(t-1)})\boldsymbol g_t^\text T(\boldsymbol\theta^{(t-1)})\mathbf S_{t-1},
\end{equation}
\end{small}
where $\alpha=1+\boldsymbol g_t^\text T(\boldsymbol\theta^{(t-1)})\mathbf S_{t-1}\boldsymbol g_t(\boldsymbol\theta^{(t-1)})$. Note that the RGN algorithm does not take advantage of the special structure.
%In addition, a damped term can be added to the Hessian matrix, i.e., $\mathbf R_t=\mathbf H_t+\delta\mathbf I$, to obtain the RLM algorithm \cite{ngia2000efficient}.
%%The derivation is similar to the RGN algorithm ( readers can also refer to \cite{ngia2000efficient} for more details).
%It should be noted that the RGN algorithm and RLM algorithm do not take advantage of the special structure and the coupling relationship of the estimated parameters.

\section{EPI based online learning algorithm}
By considering the tight coupling relationship between the linear parameters and nonlinear parameters, the VP algorithm achieves excellent performance in many applications. However the complexity of the reduced function (\ref{E-3}) and the corresponding Jacobian matrix (\ref{jac1}) prevents the classical VP algorithm directly in real-time or time-varying systems. The equivalent form presented in Section II provides a valuable eye to design an online learning algorithm.

\subsection{Algorithm}
In this subsection, we propose a novel recursive algorithm (named as REPI) based on the equivalent relationship of VP algorithm and the RGN algorithm. The proposed algorithm mainly consists of three steps.

\subsubsection{ Eliminate the linear parameters}

Instead of calculating the generalized inverse of a matrix, we employ the RLS method to estimate the linear parameters and to project them out of the objective function. Assume that $\boldsymbol \theta^{(t-1)}=(\boldsymbol a^{(t-1)}, \boldsymbol c^{(t-1)})$ are the estimated parameters using the previous $t-1$ observed data. For the new observation $y_t$ and the innovation vector $\tilde{\boldsymbol \phi}_t(\boldsymbol a^{(t-1)})=(\phi_{1}(\boldsymbol a^{(t-1)}), \cdots,\phi_{n}(\boldsymbol a^{(t-1)}))^\text T$, the procedure of estimating (or ``eliminating'') the linear parameters can be listed as follows:
\begin{equation}\label{4-1}
  \tilde{\boldsymbol c}^{(t)}=\boldsymbol c^{(t-1)}+\tilde{\boldsymbol p}_t v_t(\boldsymbol\theta^{(t-1)}),
\end{equation}
\begin{equation}\label{4-2}
  \tilde{\boldsymbol p}_t=\frac{\mathbf K_{t-1}\tilde{\boldsymbol \phi}_t(\boldsymbol a^{(t-1)})}{1+\tilde{\boldsymbol \phi}_t(\boldsymbol a^{(t-1)})^\text T\mathbf K_{t-1}\tilde{\boldsymbol \phi}_t(\boldsymbol a^{(t-1)})},
\end{equation}
where $\mathbf K_{t-1}$ is the corresponding covariance matrix of linear parameters $\boldsymbol c$.
\subsubsection{Update the nonlinear parameters according to the underlying principle of EPI}

%With the obtained linear parameters $\tilde{\boldsymbol c}^{(t)}$, we can calculate the gradient of $ v_t$ with respect to the nonlinear parameters
The gradient of $v_t$ with respect to the parameters is
\begin{small}
\begin{equation}\label{4-4}
  \frac{\partial v_t}{\partial \boldsymbol a}=[-\sum\limits_{j=1}^nc_j\frac{\partial \phi_{j}}{\partial a_1},-\sum\limits_{j=1}^nc_j\frac{\partial \phi_{j}}{\partial a_2},\cdots,-\sum\limits_{j=1}^nc_j\frac{\partial \phi_{j}}{\partial a_k}]^\text T,
\end{equation}
\end{small}
\begin{small}
\begin{equation}\label{4-41}
  \frac{\partial v_t}{\partial \boldsymbol c}=\left[-\phi_1,-\phi_2,\cdots,-\phi_n\right]^\text T.
\end{equation}
\end{small}
According to (\ref{3-10}) and (\ref{3-11}), the RGN algorithm utilizes the gradient $\boldsymbol g_t=[\frac{\partial v_t}{\partial \boldsymbol a},\frac{\partial v_t}{\partial \boldsymbol c}]$ to optimize the parameters. In such case, the coupling relationship between the linear parameters and nonlinear parameters is ignored.
\begin{figure*}[htbp]
  \centering
  \includegraphics[width=0.75\textwidth]{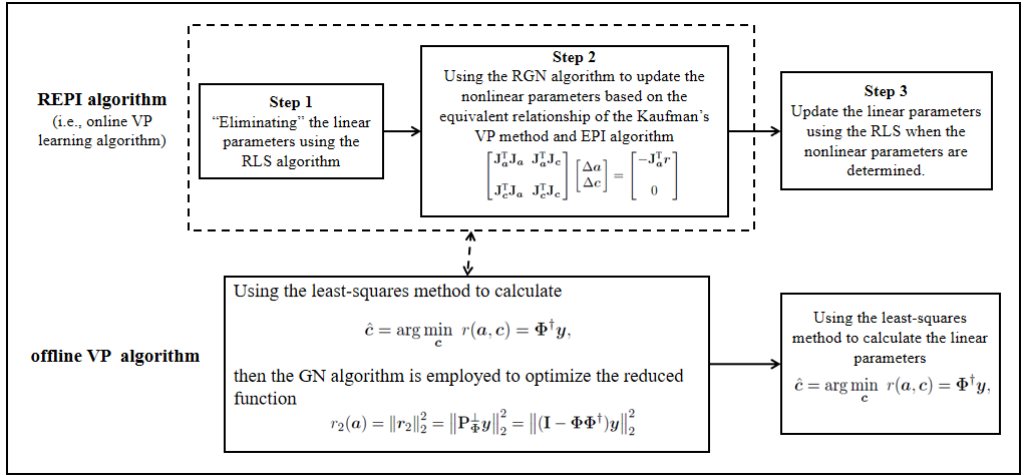}\\
  \caption{The comparison of the online VP learning algorithm (REPI) and offline VP algorithm}\label{Fig1}
\end{figure*}

Motivated by the equivalence relationship of VP algorithm and EPI algorithm, we adopt EPI principle to update the nonlinear parameters, which considers the relationship between linear parameters and nonlinear parameters. According to (\ref{E-6}), the derivation of the residual function with respect to linear parameters is equal to 0. With the linear parameters $\tilde{\boldsymbol c}^{(t)}$ calculated in step 1), the gradient with respect to nonlinear parameters in (\ref{4-4}) can be extended to
\begin{equation}\label{4-5}
  \tilde{\boldsymbol g}_t=[-\sum\limits_{j=1}^n\tilde{c}_j^{(t)}\frac{\partial \phi_{j}}{\partial a_1},\cdots,-\sum\limits_{j=1}^n\tilde{c}_j^{(t)}\frac{\partial \phi_{j}}{\partial a_k},0,\cdots,0]^\text T\in R^{k+n}.
\end{equation}

It can be seen from (\ref{E-5}) and (\ref{E-6}) that, except for the difference on the right side of the update equations, the implementation of the RGN algorithm and EPI algorithm are the same. Using the underlying principle of EPI algorithm, plugging the extended gradient (\ref{4-5}) into (\ref{3-10}) yields the update procedure of nonlinear parameters:

\begin{small}
\begin{equation}\label{4-6}
  \boldsymbol\theta^{(t)}=\boldsymbol\theta^{(t-1)}-\mathbf S_t\tilde{\boldsymbol g}_t(\boldsymbol a^{(t-1)},\tilde{\boldsymbol c}^{(t)})v_t(\boldsymbol a^{(t-1)},\tilde{\boldsymbol c}^{(t)}),
\end{equation}
\begin{eqnarray}\label{4-7}
 \mathbf S_t=\mathbf S_{t-1}-\frac{1}{\alpha_t}\mathbf S_{t-1}\boldsymbol g_t(\boldsymbol a^{(t-1)},\tilde{\boldsymbol c}^{(t)})\boldsymbol g_t(\boldsymbol a^{(t-1)},\tilde{\boldsymbol c}^{(t)})^\text T\mathbf S_{t-1},
\end{eqnarray}
\end{small}
where \begin{small}$\alpha_t=1+\boldsymbol g_t(\boldsymbol a^{(t-1)},\tilde{\boldsymbol c}^{(t)})^\text T\mathbf S_{t-1}\boldsymbol g_t(\boldsymbol a^{(t-1)},\tilde{\boldsymbol c}^{(t)})$\end{small}, $\mathbf S_t$ is the covariance matrix of estimated parameters.  Unlike the RGN method, the update direction of the linear parameters included in $\boldsymbol\theta^{(t)}$ is bypassed in this step, and it is updated using RLS method in the next step.

\subsubsection{ Update the linear parameters}

The linear parameters are updated using the RLS algorithm after the nonlinear parameters are determined at time $t$. Using the observation $y_t$ and the newest innovation vector $\boldsymbol\phi_t(\boldsymbol a^{(t)})=(\phi_{1}(\boldsymbol a^{(t)}),\cdots,\phi_{n}(\boldsymbol a^{(t)})^\text T$, the linear parameters are updated as follows:
\begin{equation}\label{4-8}
  \boldsymbol c^{(t)}=\boldsymbol c^{(t-1)}+\boldsymbol p_t v_t(\boldsymbol a^{(t)},\boldsymbol c^{(t-1)}),
\end{equation}
\begin{equation}\label{4-9}
  \boldsymbol p_t=\frac{\mathbf K_{t-1}\boldsymbol \phi_t(\boldsymbol a^{(t)})}{1+\boldsymbol \phi_t(\boldsymbol a^{(t)})^\text T\mathbf K_{t-1}\boldsymbol \phi_t(\boldsymbol a^{(t)})},
\end{equation}
\begin{equation}\label{4-10}
  \mathbf K_t=\left(\mathbf I-\boldsymbol p_t\boldsymbol \phi_t(\boldsymbol a^{(t)})^\text T\right)\mathbf K_{t-1}.
\end{equation}

\textbf{\emph{Remark:}} The idea of introducing step 1) is based on the underlying idea of the offline VP algorithm that projects out the linear parameters, resulting in a reduced function. When updating the nonlinear parameters,  the proposed REPI algorithm  considers the coupling relationship between the linear parameters and nonlinear parameters according to the equivalent form of the VP algorithm (presented in Section II). Fig.\ref{Fig1} outlines the comparison of the REPI algorithm and offline VP method.
\subsection{Convergence of the proposed REPI algorithm}
In this subsection, we give the convergence analysis of the proposed recursive algorithm. Assume that the covariance of the linear and nonlinear parameters are bounded, that is, there exist positive constants $\kappa$, $\tau$,
\begin{equation*}
s.t.,\   \left\|\mathbf S_t\right\|_2\leq\kappa, \left\|\mathbf K_t\right\|_2\leq\tau.
\end{equation*}
The following two theorems indicate that the estimation error obtained by the proposed algorithm is mean-square bounded.

\textbf{Definition 1.} \emph{An information vector $\boldsymbol\varphi(t)$ is persistently exciting \cite{johnstone1982exponential}, if there exist positive constants $\beta$, $\gamma$, and integer number $N\geq\text{dim}(\boldsymbol\varphi(t))$ such that}
\begin{small}
\begin{equation*}
  \mathbf 0<\beta\mathbf I\leq\frac{1}{N}\sum\limits_{i=0}^{N-1}\boldsymbol\varphi(t+i)\varphi(t+i)^\text T\leq\gamma\mathbf I<\infty, t>0, a.s.
\end{equation*}
\end{small}

For convenience, we give some instruction of some variables appearing in the theorems. $\boldsymbol a$ and $\boldsymbol c$ are assumed to be the true values of the parameters in the nonlinear regression models. $\boldsymbol g_i=\frac{\text d v_i}{\text d \boldsymbol \theta}(\theta^{(t-1)})$ is the gradient of residual function $v_i$ with respect to the parameters $\boldsymbol\theta=[\boldsymbol a,\boldsymbol c]\in R^{k+n}$, and $\tilde{\boldsymbol g}_i$ is the extended gradient that is calculated according to (\ref{4-5}).

\begin{thm}
Assume that $\{\boldsymbol g_t\}$ is persistently exciting, $\{\varepsilon_t,\mathcal{F}_t\}$ and $\{\delta y_t,\mathcal{F}_t\}$ are  martingale difference sequences \cite{gan2020recursive,johnstone1982exponential}, where$\{\varepsilon_t\}$ is the observation noise, $\delta y_t=\boldsymbol\phi_t(\boldsymbol a^{(t-1)})^\text T\hat{\boldsymbol c}^{(t)}-\boldsymbol\phi_t(\boldsymbol a)^\text T\boldsymbol c$, and $\mathcal{F}_t$ is a $\sigma$-algebra that is generated from the observation data up to time $t$. If $\{\varepsilon_t\}$ and $\{\delta y_t\}$ satisfy
\vspace{0.2cm}
\begin{enumerate}
  \item $E\big[\varepsilon_t|\mathcal{F}_t\big]=0,$ $E\big[\delta y_t|\mathcal{F}_t\big]=0,$ a.s.
  \vspace{0.2cm}
  \item $E\big[\varepsilon_t^2|\mathcal{F}_t\big]\leq \sigma_\varepsilon^2<\infty,$ $E\big[\delta y_t^2(t)|\mathcal{F}_t\big]\leq \sigma_y^2<\infty,$ a.s.
\end{enumerate}
\vspace{0.2cm}
then the estimation error of the nonlinear parameters obtained by the proposed recursive algorithm is mean-square bounded, i.e.,
$E\big[\|\boldsymbol a^{(t)}-\boldsymbol a\|^2\big]\leq C_1<\infty$.
\end{thm}
\begin{proof}
By the update equation (\ref{4-6}), the estimation error can be expressed as
\begin{equation}\label{4-11}
  \delta\boldsymbol\theta^{(t)}=\delta\boldsymbol\theta^{(t-1)}-\mathbf S_t\tilde{\boldsymbol g}_tv_t=\delta\boldsymbol\theta^{(t-1)}-\mathbf S_t\tilde{\boldsymbol g}_t(\varepsilon_t-\delta y_t).
\end{equation}
Let
\begin{equation}\label{4-12}
  \text T_1(t)=(\delta\theta^{(t)})^\text T\mathbf S_t^{-1}\delta\theta^{(t)},
\end{equation}
and replace $\delta\theta^{(t)}$ with (\ref{4-11}), we have
\begin{small}
\begin{equation*}
  \text T_1(t)=[\delta\theta^{(t-1)}-\mathbf S_t\tilde{\boldsymbol g}_t(\varepsilon_t-\delta y_t)]^\text T\mathbf S_t^{-1}[\delta\theta^{(t-1)}-\mathbf S_t\tilde{\boldsymbol g}_t(\varepsilon_t-\delta y_t)]
\end{equation*}
\end{small}
According to (\ref{3-8}), the recurrence of Hessian matrix is $\mathbf S_t^{-1}=\mathbf S_{t-1}^{-1}+\boldsymbol g_t \boldsymbol g_t^\text T$, then
\begin{small}
\begin{eqnarray}\label{4-13}
  &&\text T_1(t) =\text T_1(t-1)+(\boldsymbol g_t^\text T\delta\theta^{(t-1)})^2 \nonumber\\
  &&+\tilde{\boldsymbol g}_t^\text T\mathbf S_t\tilde{\boldsymbol g}_t(\varepsilon_t-\delta y_t)^2-2\tilde{\boldsymbol g}_t^\text T\delta\theta^{(t-1)}(\varepsilon_t-\delta y_t) \\
  &&\leq \text T_1(t-1)+\tilde{\boldsymbol g}_t^\text T\mathbf S_t\tilde{\boldsymbol g}_t(\varepsilon_t-\delta y_t)^2-2\tilde{\boldsymbol g}_t^\text T\delta\theta^{(t-1)}(\varepsilon_t-\delta y_t) \nonumber
\end{eqnarray}
\end{small}

Since $\{\boldsymbol g_t\}$ is persistently exciting, there exist $\beta_1,\gamma_1>0$, and an integer ($N_1\geq k+n$) such that
\begin{equation*}
  \beta_1\mathbf I\leq\frac{1}{N_1}\sum\limits_{i=0}^{N_1-1}\boldsymbol g_{t+i}\boldsymbol g_{t+i}^\text T\leq\gamma_1\mathbf I.
\end{equation*}
Therefore,
\begin{equation*}
  \beta_1N_1(k+n)\leq\sum\limits_{i=0}^{N_1-1}\left\|\boldsymbol g_{t+i}\right\|^2\leq\gamma_1N_1(k+n),
\end{equation*}
then \begin{small}$\left\|\tilde{\boldsymbol g}_{t}\right\|^2\leq\left\|\boldsymbol g_{t}\right\|^2\leq\sum\limits_{i=0}^{N_1-1}\left\|\boldsymbol g_{t+i}\right\|^2\leq\gamma_1N_1(k+n)$.\end{small} According to the bounded assumption of $\mathbf S_t$, we have
\begin{equation*}
  \tilde{\boldsymbol g}_{t}^\text T\mathbf S_t\tilde{\boldsymbol g}_{t}\leq\kappa\gamma_1N_1(k+n).
\end{equation*}

Taking the conditional expectation of both side of (\ref{4-13}) with respect to $\mathcal{F}_{t-1}$, and using the conditions 1) and 2), we have
\begin{equation*}
  E[\text T_1(t)|\mathcal{F}_{t-1}]\leq\text T_1(t-1)+\kappa\gamma_1N_1(k+n)(\sigma_\varepsilon^2+\sigma_y^2).
\end{equation*}
Since $\{\varepsilon_t\}$ and $\{\delta y_t\}$ are martingale difference sequences, we have $\{\text T_1(t)\}$ is martingale sequences \cite{gan2020recursive,johnstone1982exponential}. By the martingale convergence theorem \cite{solo1980convergence}, $\{\text T_1(t)\}$ converges to a finite random variable $\text T_0$ almost surely.

From (\ref{4-12}), we can obtain
\begin{small}
\begin{equation*}
  \text T_1(t)=(\delta\boldsymbol\theta^{(t)})^\text T\mathbf S_t^{-1}\delta\theta^{(t)}\geq||\delta\boldsymbol\theta^{(t)}||^2\sigma_\text{min}(\mathbf S_t^{-1})=\frac{1}{\kappa}||\delta\boldsymbol\theta^{(t)}||^2,
\end{equation*}
\end{small}
therefore,
\begin{equation*}
  E[||\delta\boldsymbol a^{(t)}||^2]\leq E[||\delta\boldsymbol \theta^{(t)}||^2]\leq \kappa\text T_1(t)\rightarrow\kappa\text T_0<\infty.
\end{equation*}
\end{proof}

\begin{thm}
If the innovation vector $\{\boldsymbol\phi_t(\boldsymbol a^{(t)})\}$ is persistently exciting, then the estimation error of the linear parameters obtained by the proposed algorithm is mean-square bounded.
\end{thm}
\begin{proof}
Denote $\delta \boldsymbol c^{(t)}=\boldsymbol c^{(t)}-\boldsymbol c$, and choose
\begin{equation}\label{4-14}
  \text T_2(t)=(\delta \boldsymbol c^{(t)})^\text T\mathbf K_{t-1}^{-1}\delta \boldsymbol c^{(t)},
\end{equation}
we have
\begin{eqnarray*}
  \text T_2(t)-\text T_2(t-1)=-(\delta\boldsymbol c^{(t)})^\text T(\boldsymbol\phi_t\boldsymbol\phi_t^\text T)\delta\boldsymbol c^{(t)}\leq 0.
\end{eqnarray*}
The equation holds since the recursive formula (\ref{4-10}). Thus,
\begin{equation*}
  \text T_2(t)\leq \text T_2(t-1)\leq\cdots\leq\text T_2(0).
\end{equation*}
Taking expectation of both side of (\ref{4-14}) and noting that $\mathbf K_{t-1}$ is bounded, we have
\begin{equation}\label{4-15}
  E[||\delta\boldsymbol c^{(t)}||^2]\leq\frac{\text T_2(t)}{\sigma_{\text {min}}(\mathbf K_{t-1}^{-1})}\leq\tau\text T_2(0)\leq \infty,
\end{equation}
therefore, the theorem follows.
\end{proof}
\section{Numerical illustration}
In this section, three numerical examples, including parameter estimation of a complex exponential model and fitting real-world time series using RBF-AR(X) models,  are used to verify the performance of the proposed recursive algorithm. %The experiments are carried out using MATLAB 2016b on a 2.2GHz laptop PC.
\subsection{Parameter estimation of complex exponential model}

We consider the complex exponential model with three basis functions here:
\begin{eqnarray}\label{5-1}
  y(\boldsymbol x; \boldsymbol a, \boldsymbol c)&=&c_1e^{-a_2x_1^2}\cos(a_3x_1)+c_2e^{-a_1x_1^2}\cos(a_2x_2)\nonumber\\
  &&+c_3e^{-a_4x_1^2}\sin(a_1x_3)+\epsilon,
\end{eqnarray}
where $\boldsymbol c=(c_1, c_2, c_3)^\text T$ and $\boldsymbol a=(a_1, a_2, a_3, a_4)^\text T$ are the linear and nonlinear parameters to be estimated respectively. As discussed in \cite{gan2020recursive}, the selected true parameters are $\boldsymbol c=(2, 3, 2)^\text T$ and $\boldsymbol a=(1, 1.5, 3, 0.8)^\text T$. $\epsilon$ is the Gaussian white noise with 0 mean and standard deviation 0.2, and the  variable $\boldsymbol x$ follows the standard normal distribution. Using these parameters, we randomly generate 1000 data samples. %according to the model (\ref{5-1}).

The structure of (\ref{5-1}) is assumed to be known, and we use the generated (observed) data $\{y_i\}$ to identify the parameters of the model. Four different recursive learning algorithms, including the proposed algorithm, RGN algorithm, HRGN algorithm and RVP algorithm, are employed to estimate the parameters of model (\ref{5-1}). The HRGN algorithm is a hybrid RGN algorithm, which divides the parameters into a linear part and nonlinear part and optimizes them alternately. The RVP algorithm was proposed in \cite{gan2020recursive}, which just uses the previous $p$ data points to eliminate the linear parameters. The relative errors between the estimated parameters and true parameters are used to evaluate the performances of different algorithms,
\begin{equation}\label{5-2}
  \delta\boldsymbol\theta^{(t)}=||\boldsymbol\theta^{(t)}-\boldsymbol\theta||_2/||\boldsymbol\theta||_2.
\end{equation}
\begin{figure}
  \centering
  \includegraphics[width=0.46\textwidth]{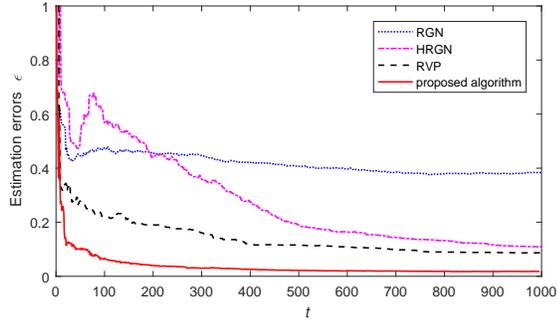}\\
  \caption{Convergence process using different recursive algorithms}\label{Fig5}
\end{figure}

For fair comparison, the observations and initial values of estimated parameters are kept the same for each algorithm. Fig. \ref{Fig5} shows the convergence process of different algorithms, and Table I lists part of the detailed results of parameter estimation during the recursive process. From Fig. \ref{Fig5} and Table I, we can observe that 1) the proposed recursive algorithm requires fewer observations to converge and achieves smallest estimation error, outperforming the other three algorithms; 2) similar to the off-line algorithm, the RGN algorithm that neglects the separable structure, is difficult to converge.

To further evaluate the efficiency and robustness of the proposed algorithm, we randomly generate 300 sets of initial values of the estimated parameters and observations. The generated initial values obey the uniform distribution, i.e.,
%$$a_1\sim\text U(0.5, 1.5),a_2\sim\text U(1, 2),a_3\sim\text U(2, 4),a_4\sim\text U(0.3, 1.3)$$
\begin{small}
\begin{eqnarray*}
  &a_1\sim\text U(0.5, 1.5),a_2\sim\text U(1, 2),a_3\sim\text U(2, 4),a_4\sim\text U(0.3, 1.3)&\\
  &c_1\sim\text U(0, 4),c_2\sim\text U(1, 5),c_3\sim\text U(0,4),&
\end{eqnarray*}
\end{small}
where $\text U(\alpha,\beta)$ represents a uniform distribution over the interval $[\alpha, \beta]$.

\begin{figure}
  \centering
  \includegraphics[width=0.46\textwidth]{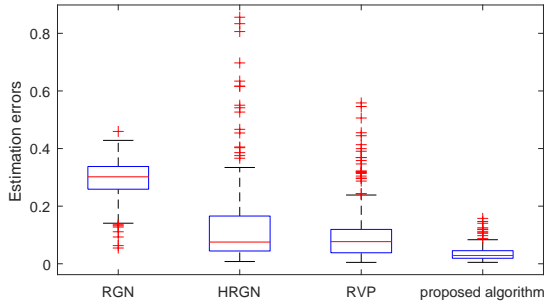}\\
  \caption{Comparison of parameter estimation errors using four different recursive algorithms for 300 runs}\label{Fig6}
\end{figure}

For each run, different algorithms start from the same initial points and the observations are kept the same. Fig. \ref{Fig6} shows the boxplot of parameter estimation errors of different recursive algorithms for 300 runs. As shown in Fig. \ref{Fig6}, we can observe that the proposed algorithm achieves the smallest average estimation error and smallest standard derivation. Both the HRGN algorithm and RVP algorithm have some large estimation errors. This mainly because that the HRGN algorithm regards linear parameters and nonlinear parameters as independent and the RVP algorithm only uses a small amount of data (here $p=10$) to eliminate the linear parameters, which makes them sensitive to the initial values of estimated parameters.

The above comparison indicates that by introducing the EPI step, the proposed recursive algorithm is more robust and efficient than the algorithms that neglect the separable structure and the algorithms that regard the linear/nonlinear parameters as independent.
\begin{table}[tp]
\centering
\begin{center}
\caption{Estimations and errors of the different recursive algorithms}
\end{center}
\label{ratiofordifferentsize}
\tabcolsep 3pt
\begin{tabular}{l|lccccccc|c }
\hline%\toprule
  &$t$ &$c_1$ &$c_2$&$c_3$ &$a_1$ &$a_2$&$a_3$&$a_4$ &$\delta(\%)$\\
\hline%\midrule\hline
\multirow{4}{*}{RGN}
&  100 &   3.36 &   1.59 &   1.71 &  1.61 &   2.418 &   1.724 &   0.74 &  47.55\\
&  200 &   3.36 &   1.69 &   1.62 &  1.42 &   2.316 &   1.724 &   0.61 &  45.59\\
&  500 &   3.27 &   1.81 &   1.62 &  1.29 &   2.140 &   1.913 &   0.49 &  40.64\\
&  1000 &  3.24 &   1.87 &   1.63 &  1.17 &   2.050 &   1.993 &   0.44 &  38.42\\
\hline
\multirow{4}{*}{HRGN}
&  100 &   2.59 &   1.49 &   1.60 &   2.11 &   2.49 &   2.65 &   2.93 &  56.81\\
&  200 &   2.86 &   2.26 &   1.47 &   1.74 &   1.89 &   2.56 &   2.68 &  44.71\\
&  500 &   2.55 &   2.60 &   1.74 &   1.34 &   1.67 &   2.60 &   1.30 &  18.96\\
&  1000 &  2.36 &   2.77 &   1.86 &   1.15 &   1.61 &   2.68 &   0.81 &  10.69\\
\hline
\multirow{5}{*}{RVP}
&  100 &   1.63 &   3.57 &   2.01 &   1.38 &   1.36 &   3.71 &   1.35 &  21.91\\
&  200 &   1.59 &   3.52 &   1.95 &   1.24 &   1.38 &   3.66 &   1.18 &  19.04\\
&  500 &   1.72 &   3.38 &   1.98 &   1.12 &   1.42 &   3.37 &   0.93 &  11.56\\
&  1000 &  1.79 &   3.28 &   1.98 &   1.08 &   1.43 &   3.29 &   0.88 &  8.67\\
\hline
&  100 &   2.20 &   3.05 &   1.81 &   1.02 &   1.58 &   3.09 &   0.63 &  6.43\\
proposed&  200 &   2.14 &   3.04 &   1.89 &   1.04 &   1.56 &   3.03 &   0.72 &  3.88\\
algorithm&  500 &   2.07 &   3.00 &   1.92 &   1.02 &   1.53 &   3.01 &   0.76 &  2.18\\
&  1000 &  2.06 &   2.99 &   1.93 &   1.01 &   1.51 &   3.01 &   0.78 &  1.79\\
\hline
True value    & &2.00 &3.00&2.00 &1.00 &1.50&3.000&0.800 &\\
\hline%\bottomrule
\end{tabular}
\end{table}
\subsection{Fitting real-world time series using RBF-AR(X) model}
The RBF-AR(X) model is a powerful statistical tool for time series analysis and system modeling, which can be expressed as
\begin{small}
\begin{displaymath}
 \left\{ \begin{array}{l}
y_t=\phi_0(\boldsymbol x_{t-1})+\sum\limits_{i=1}^p\phi_i(\boldsymbol x_{t-1})y_{t-i}+\sum\limits_{j=1}^q\varphi_j(\boldsymbol x_{t-1})u_{t-j}+e_t \\
\phi_i(\boldsymbol x_{t-1})=c_{i,0}+\sum\limits_{k=1}^mc_{i,k}\exp\{-\lambda_k\left\|\boldsymbol x_{t-1}-\boldsymbol z_k\right\|^2\}\\
\varphi_j(\boldsymbol x_{t-1})=w_{j,0}+\sum\limits_{k=1}^mw_{j,k}\exp\{-\lambda_k\left\|\boldsymbol x_{t-1}-\boldsymbol z_k\right\|^2\}\\
\boldsymbol z_k=(z_{k,1},z_{k,2},\cdots,z_{k,d})\\
%\boldsymbol y_{t-1}=(y_{t-1},y_{t-2},\cdots,y_{t-d})
\end{array} \right.
\end{displaymath}
\end{small}
where $p, q, m, d$ are the model order, $\{\lambda_k, \boldsymbol z_k, k=1,\cdots,m\}$ and $\{c_{ik}, w_{jk}, i=0,\cdots, p,j=1,\cdots, q, k=1,\cdots,q\}$ are the estimated parameters. The model is usually abbreviated as $\text{RBF-ARX}(p,q,m,d)$, when $q=0$, it degenerates to a model without exogenous inputs (denoted as $\text{RBF-AR}(p,m,d)$). The parameter identification of RBF-AR(X) model is a typical nonlinear regression problem.

In this subsection, we use the RBF-AR(X) to model two real-world time series and assume that the data samples arrive subsequently (i.e., real-time observations). Five different recursive algorithms, including the RGN, HRLM, RVP, stochastic gradient descent (SGD) method and the proposed recursive algorithm are employed to identify the model.

\subsubsection{\textbf{Glass tube drawing process}}
\begin{figure}
  \centering
  \includegraphics[width=0.35\textwidth]{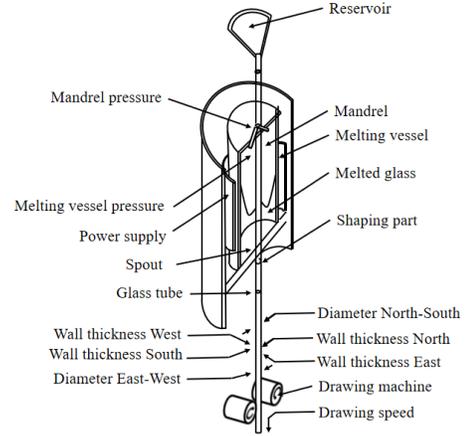}\\
  \caption{Glass tube drawing process}\label{gt}
\end{figure}
The process is outlined in Fig. \ref{gt}, and readers can refer to \cite{gan2015variable,zhu2001multivariable} for more details. Wall thickness is one of the most important quantities to be regulated, which is directly and easily affected by mandrel gas pressure and drawing speed. We used RBF-ARX(6,5,1,2) to characterize the dynamics of the process
\begin{small}
\begin{eqnarray*}
  &&y_t=\phi_0(\boldsymbol x_{t-1})+\sum\limits_{i=1}^p\phi_i(\boldsymbol x_{t-1})y_{t-i}\\
  &&+\sum\limits_{j=1}^q\left[\varphi_j^1(\boldsymbol x_{t-1}) \varphi_j^2(\boldsymbol x_{t-1})\right]\left[
                                                                                                        \begin{array}{c}
                                                                                                          u_{t-j-dl}^1 \\
                                                                                                          u_{t-j-dl}^2 \\
                                                                                                        \end{array}
                                                                                                      \right]+e_t
\end{eqnarray*}
\end{small}
where $dl$ is the time delay of the process; $\{y_t\}$, $\{u_t^1\}$ and $\{u_t^2\}$ represent the wall thickness, gas pressure and  drawing speed, respectively; the state vector $\boldsymbol x_{t-1}$ consists of the output. The data set used here (including 1269 data samples) is provided in \cite{zhu2001multivariable}. Different algorithms are applied to identified the RBF-ARX model from the first 600 subsequently observed data samples, remaining the rest to check the predictive ability of the estimated model.

\begin{figure}
  \centering
  \includegraphics[width=0.46\textwidth]{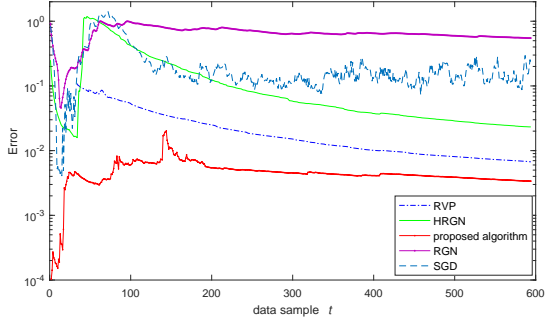}\\
  \caption{Comparison of convergence process of different recursive algorithms (using RBF-ARX to model glass tube drawing process). The error function used  here $\sum\limits_{i=1}^t(\hat{y}_i-y_i)^2/t$.}\label{Fig3}
\end{figure}

Fig. \ref{Fig3} shows the comparison of fitting results using the observed $t$ data samples. As shown in Fig. \ref{Fig3}, we can observe that the proposed recursive algorithm that considers the coupling of the linear and nonlinear parameters converges quickly, requiring fewer observations.
The RGN algorithm and SGD algorithm ignore the separable structure of RBF-AR model, which yields slower convergence rate and larger fitting error. The model obtained by the SGD algorithm is unstable and is easily affected by the new observations. The RVP algorithm eliminates linear parameters using least squares method, so it requires at least $n$ ($n$ is the number of linear parameters) observations to update the parameters of the model, and the parameters are not updated at the beginning stage. The RVP algorithm just uses a few data to tackle the coupling between parameters, which leads to the instability and slow convergence of the algorithm.

%\begin{figure}
%  \centering
%  % Requires \usepackage{graphicx}
%  \includegraphics[width=0.46\textwidth]{RBF_testfitting_1.eps}\\
%  \caption{Comparison of prediction performance of the models obtained by different methods}\label{Fig4}
%\end{figure}
\begin{table}[h]
\centering
\begin{center}
\caption{Generalization of models estimated by different algorithms (glass tube drawing process)}
\end{center}
\label{RBF_Jenkins}
\tabcolsep 3pt
\begin{tabular}{l c c c c c}
\hline%\toprule
            &RGN          & HRGN & RVP  &SGD& proposed algorithm \\
\hline%\midrule
prediction accuracy &0.5739& 0.0269&0.0149&0.0198&0.0113\\

%Ozone data&1.7608& 0.1809&0.2235&0.1410\\
\hline
\end{tabular}
\end{table}

The comparison of predictive result of models identified by the different methods is listed in Table II. From Table II, we can  observe that the model obtained by the REPI algorithm achieves best prediction performance among different recursive algorithms.

These results are similar to those of off-line algorithms for nonlinear regression problems \cite{ruhe1980algorithms, golub2003separable,gan2017some,chen2021insights}.
The comparison results show that similar to the off-line VP algorithm, the proposed algorithm is more suitable for the identification of nonlinear regression models by introducing the EPI step.
\subsubsection{\textbf{Fitting thickness ozone column data}}
\begin{figure}
  \centering
  \includegraphics[width=0.48\textwidth]{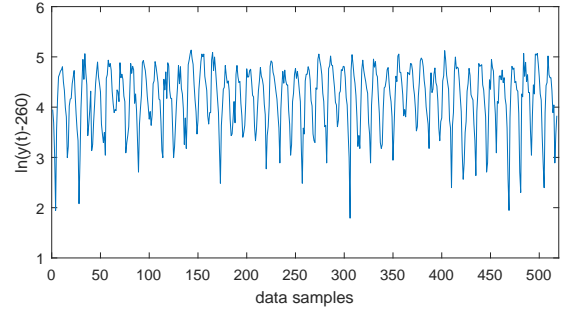}\\
  \caption{The transformed thickness ozone column data}\label{ozone}
\end{figure}

The thickness ozone column data collects 518 records of the mean thickness ozone column in Arosa, Switzerland. We do the transformation ($\tilde{y}(t)=\ln(y(t)-260)$) as \cite{chen2018ins} to make it more nearly symmetric and to stabilize the variance.
%\begin{equation*}
%  \tilde{y}(t)=\ln(y(t)-260).
%\end{equation*}
The transformed data are shown in Fig. \ref{ozone}. It is also assumed that the data arrive subsequently. We use the RBF-AR(5,1,2) to fit the data, and employ four different recursive algorithms to identify the model. The previous 300 data points are used to identify the model, remaining the rest data to check the  predictive ability of estimated models.
\begin{figure}
  \centering
  \includegraphics[width=0.48\textwidth]{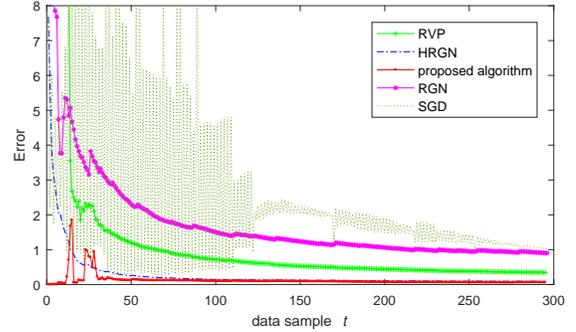}\\
  \caption{Comparison of convergence process of different recursive algorithms using RBF-AR model to fit thickness ozone column data.}\label{Fig7}
\end{figure}

Fig. \ref{Fig7} shows the convergence process of different algorithms. As shown in Fig. \ref{Fig7}, the proposed algorithm requires fewer data samples to converge than the other four algorithms. The HRGN algorithm and RVP algorithm that take advantage of the separable structure of RBF-AR model perform better than the RGN algorithm and SGD algorithm. However, the HRGN algorithm ignores the relationship between linear parameters and nonlinear parameters and the RVP algorithm only uses a small amount of data samples to eliminate linear parameters. These limit their learning performances. The SGD algorithm converges slowly and fluctuates greatly for new observations. The detailed comparisons of prediction result of different models estimated by the five recursive algorithms are listed in Table III.

\begin{table}[h]
\centering
\begin{center}
\caption{Generalization of models identified by different algorithms (thickness ozone column data)}
\end{center}
\label{RBF_ozone}
\tabcolsep 3pt
\begin{tabular}{l c c c c c}
\hline%\toprule
            &RGN          & HRGN & RVP  &SGD& proposed algorithm \\
\hline%\midrule
%prediction accuracy &1.96& 0.1476&0.0568&0.7033&0.0450\\

prediction accuracy&1.7608& 0.1809&0.2235&2.0837&0.1410\\
\hline
\end{tabular}
\end{table}

The above comparison results indicate the high efficiency of the proposed recursive algorithm. By introducing the EPI step, the proposed algorithm can adjust the parameters of the model appropriately by handling the coupling relationship between linear parameters and nonlinear parameters.

\section{Conclusion}
Nonlinear regression models have been widely used in scientific research and engineering applications. The parameter identification, especially for online system, is quite challenging because of the existence of nonlinear parameters. Empirical results have proven that the VP algorithm is valuable in solving such problems. However, previous research on the VP algorithm is limited to offline problems. For online identification of nonlinear regression problems, researchers often ignore the separable structure.

In this paper, we study the recursive algorithms for nonlinear regression models and propose an efficient recursive algorithm by analyzing the equivalent form of VP algorithm. Similar to the offline VP algorithm, the proposed algorithm is more efficient and robust than the previous RGN algorithm and HRGN algorithm by introducing the EPI step. Moreover, we prove that the proposed algorithm is mean-square bounded. This study is of great important for the development of real-time systems  because of the widespread of nonlinear regression models in the real-world applications.

\bibliographystyle{IEEEtran}
\bibliography{IEEEexample,bibfile}

\end{document}